\documentclass{birkjour}
\usepackage{amsmath}
\usepackage{amssymb}
\usepackage[english]{babel}

\newtheorem{Th}{Theorem}

\newtheorem{Prop}{Proposition}

\newtheorem*{Rem*}{Remark}

\newcommand{\N}{\mathbb{N}}
\newcommand{\Z}{\mathbb{Z}}
\newcommand{\Q}{\mathbb{Q}}
\newcommand{\R}{\mathbb{R}}

\newcommand{\eps}{\varepsilon}

\newcommand{\Span}{{\rm Span}}

\newcommand{\calN}{{\mathcal{N}}}

\newcommand{\cale}{{\mathcal{E}}}

\newcommand{\om}{\omega}
\newcommand{\pha}{\varphi}
\newcommand{\te}{\vartheta}
\newcommand{\sig}{\sigma}
\newcommand{\vareta}{\eta}

\newcommand{\unN}{\{1,\ldots,N\}}
\newcommand{\zeroN}{\{0,\ldots,N\}}

\newcommand{\norm}[1]{\lVert#1\rVert}

\title{Nesterenko's criterion when the small linear forms oscillate}
\author{St\'ephane Fischler}
\date{\today}
\address{%
Equipe d'Arithm\'etique et de G\'eom\'etrie Alg\'ebrique, 
Universit\'e Paris-Sud, B\^atiment 425,
91405 Orsay Cedex, France}

\begin{document}



\maketitle

\begin{abstract}
In this paper we generalize Nesterenko's criterion to the case where the small linear forms have an oscillating behaviour (for instance given by the saddle point method). This criterion   provides both a lower bound for the dimension of the vector space spanned over the rationals by a family of real numbers, and a measure of simultaneous approximation to these numbers (namely, an upper bound for the irrationality exponent if 1 and only one other number are involved). As an application, we prove an explicit measure of simultaneous approximation to $\zeta(5)$,  $\zeta(7)$,  $\zeta(9)$, and  $\zeta(11)$, using Zudilin's proof that at least one of these numbers is irrational.
\end{abstract}


\section{Introduction}

To investigate on Diophantine properties of real numbers $\xi_1$, \ldots, $\xi_r$, a strategy is to construct small linear forms in 1,  $\xi_1$, \ldots, $\xi_r$ with integer coefficients. This is the only known way to study this problem if  $\xi_1$, \ldots, $\xi_r$ are, for instance, values of Riemann zeta function $\zeta(s) = \sum_{n=1}^\infty \frac1{n^s}$ at odd integers $s \geq 3$. 

In more precise terms, linear forms $\ell_{0,n} + \ell_{1,n}\xi_1 + \ldots + \ell_{r,n} \xi_r$ are constructed, with absolute value $\leq \alpha^{n+o(n)}$  as $n \to \infty$ and coefficients $\ell_{j,n}\in \Z$ such that $|\ell_{j,n}|\leq \beta^{n+o(n)}$ (where $0 < \alpha  < 1 < \beta$). Then:
\begin{itemize}
\item[$\bullet$] If $\ell_{0,n} + \ell_{1,n}\xi_1 + \ldots + \ell_{r,n} \xi_r \neq 0$ for infinitely many $n$, the subgroup $\Z + \Z \xi_1 + \ldots + \Z \xi_r$ of $\R$ is not discrete so that at least one number among $\xi_1$,\ldots, $\xi_r$ is irrational.
\item[$\bullet$] If $\alpha \beta < 1$ and $\ell_{0,n} + \ell_{1,n}\xi_1 + \ldots + \ell_{r,n} \xi_r \neq 0$ for infinitely many $n$, the subspace spanned over $\Q$ by  1, $\xi_1$,\ldots, $\xi_r$ has dimension at least~3, so that 1, $\xi_i$ and $\xi_j$ are $\Q$-linearly independent for some $i$, $j$ with $1 \leq i < j \leq r$. An elementary proof of this result can be found in \cite{SFZu} (Proposition 1, \S 2.2). 
\item[$\bullet$] If   the  linear forms are not too small, namely $\geq \alpha^{n+o(n)}$, a result of Nesterenko \cite{Nesterenkocritere} implies the following two properties (see the remark after Theorem \ref{thcrit} below for a precise formulation):
\begin{itemize}
\item[$(i)$] A lower bound $1 -  \frac{\log \alpha}{\log \beta}$   on the dimension of the vector space spanned over the rationals by 1, $\xi_1$, \ldots, $\xi_r$. This linear independence criterion is one of the main tools in the proof (\cite{RivoalCRAS}, \cite{BR}) that $\zeta(s) \not\in \Q$ for infinitely many odd integers $s \geq 3$.
\item[$(ii)$] A measure of simultaneous approximation to  $\xi_1$, \ldots, $\xi_r$; if \mbox{$r=1$,} this is an upper bound $\mu(\xi_1) \leq 1 - \frac{\log \beta}{\log \alpha}$ on the irrationality exponent of $\xi_1$. For instance, Ap\'ery proved at the same time \cite{Apery} that $\zeta(3)\not\in \Q$ and $\mu(\zeta(3)) \leq 13.41\ldots$.  
\end{itemize}
The assumption that the linear forms are not too small is very important here, and it cannot be omitted. It can be weakened:   Nesterenko proves such results for linear forms with  absolute values between $\alpha_1^{n+o(n)}$ and $\alpha^{n+o(n)}$, where $0 < \alpha_1 < \alpha < 1$, but the conclusion is then weaker too.
\end{itemize}

\bigskip

In this paper we generalize Nesterenko's results $(i)$ and $(ii)$ (without weakening the conclusion) to the case where the linear forms behave essentially like $\alpha^{n+o(n)} \cos(n \om + \pha)$ with $\om, \pha\in \R$. This is an interesting situation because the saddle point method is often applied to obtain asymptotic estimates for the linear forms, and it typically produces this kind of behaviour. Our main result is the following.

\begin{Th} \label{thcrit}
Let $r \geq 1$, $ \xi_1 , \ldots,  \xi_r \in \R$, $\alpha, \beta, \om , \pha\in \R$. Assume that \mbox{$0 < \alpha < 1$,} $\beta > 1$, and either $\om \not\equiv 0 \bmod \pi$ or $\pha \not\equiv \frac{\pi}2 \bmod \pi$. For any $n \geq 1$, let $\ell_{0,n} , \ldots,  \ell_{r,n}  \in \Z$ be such that, as $n\to \infty$:
$$\max_{0\leq i \leq r} |\ell_{i,n}| \leq \beta^{n+o(n)} $$
and
\begin{equation} \label{eqcrit}
|\ell_{0,n} + \ell_{1,n}\xi_1 + \ldots + \ell_{r,n} \xi_r| = \alpha^{n+o(n)}\Big(| \cos(n \om + \pha)|+o(1)\Big).
\end{equation}
Then:
\begin{itemize}
\item[$(i)$] We have  $\dim_\Q \Span_\Q(1,  \xi_1 , \ldots,  \xi_r)\geq 1 - \frac{\log \alpha}{\log \beta} $. 
\item[$(ii)$] For any $\kappa > 1 - \frac{\log \beta}{\log \alpha} $  and any $q,p_1,\ldots,p_r\in \Z$ with $q > 0$ sufficiently large in terms of $\kappa$, we have
$$\max\Big(\Big|\xi_1 - \frac{p_1}q\Big|,\ldots, \Big|\xi_r - \frac{p_r}q\Big|\Big) \geq \frac1{q^\kappa}.
$$
\end{itemize}
\end{Th}

When $(\om,\pha)=(0,0)$ these are exactly Nesterenko's above-mentioned results \cite{Nesterenkocritere}. If $e^{i \om}$ and $e^{i \pha}$ are algebraic numbers, a  very concise remark of Sorokin \cite{SorokinZR} (which we expand in \S \ref{subsecdemcrit}) provides another proof of Theorem \ref{thcrit}, based upon lower bounds for linear forms in logarithms.

As the proof shows, the cosine may be replaced in Theorem \ref{thcrit} (and also in Proposition \ref{prop1} below) with any continuous periodic function, which vanishes only at finitely many points within each period.

\bigskip

Theorem \ref{thcrit} can be used in the following situation (see \S \ref{subsec22} for other possible applications). Zudilin has proved \cite{Zudilinonze} that among $\zeta(5)$, $\zeta(7)$, $\zeta(9)$,  and $\zeta(11)$, at least one is irrational (refining upon Rivoal's result \cite{vingtetun}). He proceeds by constructing small linear forms in 1 and these numbers. He  applies the saddle point method  to prove the estimate \eqref{eqcrit}  (see \S \ref{subsec21} below), and deduce that the linear form is non-zero for infinitely many $n$, thereby proving the irrationality of  at least one number among $\zeta(5)$, $\zeta(7)$, $\zeta(9)$,  and $\zeta(11)$. Using  Theorem \ref{thcrit}  we obtain a quantitative version of this result:

\begin{Th} \label{thzeta}
For any $q,p_5, p_7,p_9,p_{11}\in \Z$ with $q > 0$ sufficiently large we have:
$$\max\Big(\Big|\zeta(5) - \frac{p_5}q\Big|,   \Big|\zeta(7) - \frac{p_7}q\Big|, \Big|\zeta(9) - \frac{p_9}q\Big|,   \Big|\zeta(11) - \frac{p_{11}}q\Big| \Big) \geq \frac1{q^{438.23}}.
$$
\end{Th}

In particular, in the (very unlikely) case where  $\zeta(5)$, $\zeta(7)$ and $\zeta(9)$ would be rational numbers, this implies $\mu(\zeta(11)) \leq 438.23$.

\bigskip

The proof of  Theorem \ref{thcrit} (see \S  \ref{subsecdemcrit}) relies on applying   Nesterenko's results to a subsequence $\ell_{0,\psi(n)} + \ell_{1,\psi(n)}\xi_1 + \ldots + \ell_{r,\psi(n)} \xi_r$, where $\psi $ is  given by  the following proposition (with \mbox{$N=1$}).

\begin{Prop} \label{prop1} Let $N \geq 1$, $\om_1,\ldots,\om_N, \pha_1,\ldots,\pha_N\in \R$. Assume there exist infinitely many integers $n$ such that, for any $i \in \unN$, $n \om_i + \pha_i \not\equiv \frac{\pi}2 \bmod \pi$. Then there exist $\eps, \lambda > 0$ and an increasing function $\psi : \N \to \N$ such that $\lim_{n\to \infty} \frac{\psi(n)}{n} = \lambda$ and, for any $n \in \N$ and any $i\in \unN$, $|\cos(\psi(n)\om_i + \pha_i)| \geq \eps$.
\end{Prop}

This proposition   is a consequence of Kronecker-Weyl's  equidistribution theorem (see \S \ref{subsec1}). The assumption that for infinitely many  $n$,   $n \om_i + \pha_i \not\equiv \frac{\pi}2 \bmod \pi$
 for any $i  $ is of course necessary because if $\psi$ exists, all $n \in \psi(\N)$ have this property. If $N=1$, it is equivalent to $(\om_1,\pha_1) \not\equiv (0 \bmod \pi, \frac{\pi}2 \bmod \pi)$. For $N=2$, it is equivalent to:
 \begin{equation} \label{eqNde}
 \left\{
 \begin{array}{l}
 (\om_1,\pha_1) \not\equiv (0 \bmod \pi, \frac{\pi}2 \bmod \pi)\\
 (\om_2,\pha_2) \not\equiv (0 \bmod \pi, \frac{\pi}2 \bmod \pi)\\
 (\om_1, \om_2,\pha_1,\pha_2) \not\equiv ( \frac{\pi}2  \bmod \pi, \frac{\pi}2 \bmod \pi,  \frac{\pi}2  \bmod \pi, 0 \bmod \pi)\\
 (\om_1, \om_2,\pha_1,\pha_2) \not\equiv ( \frac{\pi}2  \bmod \pi, \frac{\pi}2 \bmod \pi, 0 \bmod \pi,   \frac{\pi}2  \bmod \pi)
 \end{array}\right.
\end{equation}
 In all Diophantine applications we have in view (see \S \ref{subsec23}), $N$ is fixed (and even $N=1$ for the above-mentioned results), so that this assumption is not too difficult to check. It should be pointed out that in Proposition \ref{prop1}, no Diophantine condition is assumed on $\om_i$ or $\pha_i$: for instance, we don't assume $\pi$, $\om_1$, \ldots, $\om_N$ to be $\Q$-linearly independent. This is very useful because such a condition could be very difficult to check (whereas approximate values of these numbers are often easily computed, making it very easy to check conditions such as \eqref{eqNde}).

\bigskip



\noindent{\bf Acknowledgements:} I am very thankful to Francesco Amoroso and Michel Waldschmidt for reading carefully a preliminary version of this paper, and making helpful remarks; in this respect my gratitude goes to the referee too. I would like also to thank warmly Tanguy Rivoal for his advice, and  Wadim Zudilin for pointing out Sorokin's paper to me.

\section{Proofs}

\subsection{Proof of Proposition \ref{prop1}} \label{subsec1}

Let  $p : \R \to \R/\Z , x \mapsto x \bmod \Z $ denote the canonical surjection. For $x,y\in\R$ with $x<y$, we denote by $[x,y]\subset\R/\Z$ the image under $p$ of the segment between $x$ and $y$, that is the set of all $p(z)$ with $x\leq z \leq y$. Its normalized Haar measure is $\min(y-x,1)$.

Let $s$ be a positive integer.  For $x_1,\ldots,x_s,$ $y_1,\ldots,y_s \in \R$ such that $x_i < y_i$, let $K = \prod_{i=1}^s [x_i,y_i]  \subset (\R/\Z)^s  $.  Given  $\te_1,\ldots, \te_s \in \R$, we denote by  $\calN$  the set of all $n\in \N$ such that $(n\te_1 \bmod \Z ,\ldots,n\te_s   \bmod \Z ) \in K$. Then Kronecker-Weyl's  equidistribution  theorem (see \cite{KuipersNiederreiter}, Chapter~1.6) asserts that, if $1,\te_1,\ldots,\te_s$ are $\Q$-linearly independent:
$$\lim_{k \to  \infty} \frac{| \{n \in \calN , \, n \leq k\}|}{k} = \prod_{i=1}^s \min( y_i - x_i , 1) .$$
Denoting by   $\psi_0 (n)$   the $n$-th element of $\calN$  in increasing order (so that $\psi_0 : \N \to \calN$ is an increasing bijection), this implies
\begin{equation} \label{eqKW}
\lim_{n \to \infty} \frac{\psi_0(n)}{n} = \prod_{i=1}^s \max(  (y_i - x_i)^{-1}, 1)
\end{equation}
by letting $k  = \psi_0 (n)$ and taking reciprocals. Equation \eqref{eqKW} will be the main tool in the proof of Proposition \ref{prop1}.

\bigskip

To illustrate the ideas in a simpler case,   let us prove  Proposition \ref{prop1} first when $N=1$. If $\om_1/\pi$ is a rational number $c/d$, then it is enough to choose $\psi(n) = nd + a$ for a suitable $a \in \{1,\ldots,d\}$. If $\om_1/\pi$ is irrational, let $K = [-\frac{\pha_1}{\pi} - \frac14, -\frac{\pha_1}{\pi}  + \frac14]  \subset \R/\Z$ and $\calN$ be the  set of all $n\in \N$ such that $\frac{ n \om_1}{\pi}  \bmod \Z   \in K$. 
Denoting by   $\psi_0 (n)$   the $n$-th element of $\calN$, Eq.  \eqref{eqKW}  yields $\lim_{n \to \infty} \frac{\psi_0(n)}{n}  = 2$. Moreover for any $n$ we have 
$$\Big|\cos(\psi_0 (n) \om_1 + \pha_1)\Big| = \Big|\cos\Big(\pi\Big( \psi_0 (n) \frac{ \om_1}{\pi} + \frac{\pha_1}{\pi} \Big)\Big)\Big| \geq \frac{\sqrt 2}2.$$
 This concludes the proof when $N=1$.

\bigskip

Let us come now to the proof of Proposition \ref{prop1} for any $N$, starting with  a special case: we assume   that $\frac{\om_1}{\pi}, \ldots, \frac{\om_N}{\pi}$ are irrational numbers. Let $E$ denote the $\Q$-subspace of $\R$ spanned by $1, \frac{\om_1}{\pi}, \ldots, \frac{\om_N}{\pi}$. Since $\frac{\om_1}{\pi} \not\in \Q$, we have $\dim_\Q E \geq 2$; let $(1,\te_1,\ldots,\te_s)$ denote a basis of $E$, with $s \geq 1$. For any $i\in \unN$, we have $\frac{\om_i}{\pi} = \sum_{j=0}^s r_{i,j} \te_j$ with $r_{i,j} \in \Q$, where we let $\te_0 = 1$. Let $D $ be a positive integer such that $D r_{i,j} \in \Z$ for any $i,j$. Then we have
 \begin{equation} \label{eqdcp}
\Big|  \cos(Dn\om_i + \pha_i ) \Big|  = \Big|  \cos\Big(\pi\Big(\sum_{j=1}^s D r_{i,j} n \te_j + \frac{\pha_i}{\pi}\Big)\Big) \Big|
\end{equation}
 because $D r_{i,0}n\te_0 \in \Z$.
 
 For any $i\in \unN$, let $\Delta_i \subset (\R/\Z)^s$ be the set of all $(\sig_1,\ldots,\sig_s)$ such that $ (\sum_{j=1}^s D r_{i,j} \sig_j ) + \frac{\pha_i}{\pi} - \frac12  \in \Z$. Since $(Dr_{i,1} ,\ldots,Dr_{i,s})\in \Z^s \setminus \{(0,\ldots,0)\}$ because $\frac{\om_i}{\pi} \not\in \Q$, $\Delta_i$ is a  finite union of translated tori of dimension $s-1$, and a proper compact subset of  $(\R/\Z)^s$. There exists a point $(z_1,\ldots, z_s) \in \R^s$  and a (small) positive real number $\eta$ such that $K = \prod_{i=1}^s [z_i -\eta, z_i + \eta]  \subset ( \R/\Z )^s$ is disjoint from $\Delta_1\cup\ldots\cup \Delta_N$.  Let $\calN$ be the set of all $n\in \N$ such that $(n\te_1  \bmod \Z ,\ldots,n\te_s \bmod \Z ) \in K$, and $\psi_0 (n)$ denote the $n$-th element of $\calN$ (in increasing order). Since $1,\te_1,\ldots,\te_s$ are $\Q$-linearly independent, Kronecker-Weyl's equidistribution theorem \eqref{eqKW} yields
$$\lim_{n \to \infty} \frac{\psi_0(n)}{n} = (2 \eta)^{-s} > 0.$$
 Moreover, since $K$, $\Delta_1$,  \ldots, $  \Delta_N $ are  compact subsets there exists $\vareta' >0$ such that, for any $n \in \N$ and any $i\in \unN$, $\norm{ ( \sum_{j=1}^s D r_{i,j} \psi_0(n) \te_j  ) + \frac{\pha_i}{\pi}-\frac12 }_\Z \geq \vareta'$ (where $\norm{x}_\Z$ is the distance of $x\in\R$ to $\Z$). Using Eq. \eqref{eqdcp}, this provides $\eps>0$ such that $| \cos(D\psi_0(n)\om_i + \pha_i )| \geq \eps$, and by letting $\psi(n) = D\psi_0(n)$ this concludes the proof of  Proposition \ref{prop1} if  $\frac{\om_1}{\pi}, \ldots, \frac{\om_N}{\pi}\not\in \Q$.
 
 \bigskip
 
 Let us deduce the general case from this special case. Reordering the pairs $(\om_i,\pha_i)$ if necessary, we may assume that for some $N' \in \zeroN$ we have $\frac{\om_1}{\pi}, \ldots, \frac{\om_{N'}}{\pi}\not\in \Q$ and $\frac{\om_{N'+1}}{\pi}, \ldots, \frac{\om_N}{\pi}\in \Q$. Let $d\geq 1$ be a common denominator of $\frac{\om_{N'+1}}{\pi}, \ldots, \frac{\om_N}{\pi}$ and for any $i$, let $\cale_i$ be the set of all $k \in \N$ such that $k\om_i + \pha_i \equiv \frac{\pi}2 \bmod \pi$. Then $\cale_i$ has at most one element for $i\leq N'$, and $\cale_i$ is a union of residue classes mod $d$ for $i > N'$. By assumption $\N \setminus (\cale_1 \cup \ldots \cup \cale_N)$ is infinite, so that there exists $a\in \N$ such that for any $k$ sufficiently large with $k \equiv a \bmod d$, we have $k\not \in \cale_1 \cup \ldots \cup \cale_N$. For any $ i > N'$, the number $|\cos((nd+a)\om_i+\pha_i)|$ is positive and independent from $n$ (since $d\om_i \in \pi\Z$). If $N' =0$ this concludes the proof by letting $\psi(n) = nd+a$. Otherwise we apply the special case of  Proposition \ref{prop1}  proved above  to $ \cos((nd+a)\om_i+\pha_i) = \cos(nd\om_i+a\om_i+\pha_i) $ for $1 \leq i \leq N'$, that is with $\om'_1=d\om_1$, \ldots, $\om'_{N'}=d\om_{N'}$, $\pha'_1 = a\om_1+\pha_1$,\ldots,  $\pha'_{N'} = a\om_{N'}+\pha_{N'}$. We obtain in this way an increasing function $\psi_0$, and letting $\psi(n) = \psi_0(n)d+a$ concludes the proof of  Proposition~\ref{prop1}.

\begin{Rem*}
If $\frac{\om_1}{\pi}, \ldots, \frac{\om_N}{\pi}$ are irrational numbers, applying Kronecker-Weyl's theorem with more general subsets $K$ enables one to obtain  $\psi_0$  such that $\lim_{n \to \infty} \frac{\psi_0(n)}{n}$ is arbitrarily close to 1 (because this is the inverse of the measure of $K$). This leads to a control upon $\lambda = 
\lim_{n \to \infty} \frac{\psi (n)}{n}$ in terms of the common denominator $D$ in this case. If 1, $\frac{\om_1}{\pi}$, \ldots, $\frac{\om_N}{\pi}$ are $\Q$-linearly independent then we can take $D=1$, so that $\lambda$ can be chosen arbitrarily close to~1. However we did not try to go any further in this direction (nor to get a lower bound for $\eps$) because this is completely   useless for the applications we have in view.
\end{Rem*}

\subsection{Proof of Theorem \ref{thcrit}} \label{subsecdemcrit}

Since $\om \not\equiv 0 \bmod \pi$ or $\pha \not\equiv \frac{\pi}2 \bmod \pi$, there are   infinitely many integers $n$ such that  $n \om + \pha  \not\equiv \frac{\pi}2 \bmod \pi$. Applying Proposition \ref{prop1} (with $N=1$)  yields  $\eps, \lambda > 0$ and an increasing function $\psi : \N \to \N$ such that $ \psi(n)  = \lambda n + o(n)$  as $n \to \infty$ and $|\cos(\psi(n)\om  + \pha )| \geq \eps$  for any $n$.  Therefore Eq. \eqref{eqcrit} yields
$$|\ell_{0,\psi(n)} + \ell_{1,\psi(n)}\xi_1 + \ldots + \ell_{r,\psi(n)} \xi_r| = \alpha^{ \psi(n)(1+ o(1))} =   \alpha^{ \lambda n + o(n)}$$
and we have also $\max_{0\leq i \leq r} |\ell_{i,\psi(n)}| \leq \beta^{ \lambda n + o(n)}$. Therefore Nesterenko's results \cite{Nesterenkocritere} (that is, the special case of Theorem \ref{thcrit} where $\om = \pha = 0$) apply to the sequence $ \ell_{0,\psi(n)} + \ell_{1,\psi(n)}\xi_1 + \ldots + \ell_{r,\psi(n)} \xi_r$, with $\alpha^\lambda$ and $\beta^\lambda$ instead of $\alpha$ and $\beta$. This provides exactly the same conclusions $(i)$ and $(ii)$ because $\frac{\log \alpha^\lambda}{\log \beta^\lambda} = \frac{\log \alpha }{\log \beta }$. 

\begin{Rem*} Another proof of Theorem \ref{thcrit} under the
additional assumption that $e^{i \om}$ and $e^{i \pha}$ are
algebraic numbers  is provided by the following remark of Sorokin (p. 823
of \cite{SorokinZR}): 
$$\lim_{n \to \infty \atop n \in S} |\cos(n \om + \pha)|^{1/n} = 1$$
 in this case, if $\om \not\equiv 0 \bmod \pi$ or $\pha \not\equiv \frac{\pi}2 \bmod \pi$,
 where $S$ is the set of all $n$ such that
$\cos(n \omega + \varphi) \neq 0$. Indeed, it follows from
Gel'fond's lower bound for linear forms in logarithms (see
Theorem 4.1 of \cite{EMS}, p. 179) that for any $\eps > 0$ and any $k$,
$n$ sufficiently large in terms of $\eps$, we have
$$ | n \om + \pha - \frac{\pi}2 - k \pi | \geq (1-\eps)^n$$
provided the left hand-side is non-zero.
\end{Rem*}

\subsection{Simultaneous approximation to zeta values} \label{subsec21}

This subsection is devoted to a proof of Theorem \ref{thzeta}. We use Zudilin's linear forms (\cite{Zudilinonze}; see also \cite{Zudilincinqaout} for further details), which can be written as
\begin{equation} \label{eqzu}
S_n = \frac12 \frac{\prod_{j=1}^{10} ((13+2j)n)!}{(27n)!^6} \sum_{k=1}^{\infty}  \frac{{\rm d}^2}{{\rm d}t^2} \Big( (37n+2t) \frac{(t-27n)_{27n}^3 (t+37n+1)_{27n}^3}{\prod_{j=1}^{10} (t+(12-j)n)_{(13+2j)n+1}} \Big)_{|t = k},
\end{equation}
where the second derivative is taken at $t=k$ and $(\alpha)_p = \alpha(\alpha+1)\ldots(\alpha+p-1)$ is Pochhammer's symbol. This sum can be written as a linear form 
$$S_n = \ell_{0 ,n} +  \ell_{5,n} \zeta(5) + \ell_{7,n} \zeta(7) + \ell_{9,n} \zeta(9) + \ell_{11,n} \zeta(11)$$
with rational coefficients $\ell_{i,n}$. Zudilin deduces from  the saddle point method  that
\begin{equation} \label{eqzuasy}
S_n = e^{-C_0 n + o(n)} |\cos(n \om+ \pha) + o(1)|, 
\end{equation}
$$
\mbox{ with } \om \not\equiv 0 \bmod \pi \mbox{ and } C_0 = 227.58019641\ldots,
$$
and constructs a common denominator $D_n$ of the rational numbers  $\ell_{0 ,n}$, $ \ell_{5,n}$, $ \ell_{7,n}$, $\ell_{9,n}$, and $ \ell_{11,n}$, such that
$$D_n =  e^{ C_1 n + o(n)}  \mbox{ with } C_1 = 226.24944266\ldots$$
Then $D_n S_n$ is a linear form in  1, $\zeta(5)$, \ldots, $\zeta(11)$, with integer coefficients; as $n \to \infty$,   it tends to 0 because $\alpha = e^{C_1 - C_0} < 1$,
and is non-zero for infinitely many $n$ thanks to \eqref{eqzuasy}. This proves that among the four numbers $\zeta(5)$, \ldots, $\zeta(11)$, at least one is irrational.

To prove Theorem \ref{thzeta}, we also need the following upper bound on the coefficients $\ell_{j,n}$, which can be proved along the same lines as Lemma 4 of \cite{BR}:
$$\max_{j \in \{0, 5,7,9,11\}}  |\ell_{j,n} | \leq 2^{513 n + o(n)}, \mbox{ since } 3(27+37+27) + \sum_{j=1}^{10} (13+2j) = 513.$$
This allows us to apply Theorem \ref{thcrit} with $\beta =  e^{ C_1} 2^{513}$, so that $1 - \frac{\log \beta}{\log \alpha}   = 438.22134\ldots$. This concludes the proof of Theorem \ref{thzeta}.

\section{Other possible applications} \label{sec3}

\subsection{Applications of the linear independence criterion}  \label{subsec22}
 
 For any $s\geq 2$, let $i_s$ denote the least odd integer $i\geq 3$ such that 
 $$\dim_\Q \Span_\Q (1,\zeta(3),\zeta(5),\ldots, \zeta(i)) \geq s.$$ Since $\zeta(3),\zeta(5),\ldots$ span an infinite-dimensional $\Q$-vector space (\cite{BR}, \cite{RivoalCRAS}), $i_s$ exists for any $s$. Ap\'ery's result that $\zeta(3)\not\in \Q$ means $i_2=3$; conjecturally, $i_s=2s-1 $ for any $s$. Ball-Rivoal's construction yields an upper bound on $i_s$, which has been improved for small values of $s$, namely  $i_3 \leq 139$ and $i_4 \leq 1961$ (\cite{SFZu}, refining upon previous bounds due to Zudilin \cite{Zudilincentqc}).  

Now let $j_s$ denote  the least odd integer $j\geq 5$ such that 
$$\dim_\Q \Span_\Q (1,\zeta(5),\zeta(7),\ldots, \zeta(j)) \geq s.$$ 
The trivial remark that $j_s \leq i_{s+1}$ yields $j_2 \leq 139$ and $j_3 \leq 1961$. Now Zudilin's result \cite{Zudilinonze} is an important improvement of the former bound, namely $j_2 \leq 11$. The linear forms in 1, $\zeta(5)$, $\zeta(7)$, \ldots, $ \zeta(j)$ he constructs (in the style of \eqref{eqzu} above: see also  \cite{Zudilincentqc}) have the asymptotics \eqref{eqcrit}, with $\om \not\equiv 0 \bmod \pi$ in general. Theorem \ref{thcrit} enables one to deduce from it an upper bound on $j_s$, which should be better (for a fixed value of $s $) than the one derived from the trivial bound $ j_s \leq i_{s+1}$. In particular the bound $j_3 \leq 1961$ can probably be improved in this way (note however that for proving that at least three numbers in a family are $\Q$-linearly independent, Theorem \ref{thcrit} can be replaced with 
 Proposition 1 of \cite{SFZu}, as explained at the beginning of the introduction: a lower bound for the linear forms is not necessary in this case).

\subsection{Further generalizations of Nesterenko's criterion}  \label{subsec23}

Nesterenko's results $(i)$ and $(ii)$ (with $\om = \pha = 0$ in Eq. \eqref{eqcrit})  have been generalized in several directions, namely:
\begin{itemize}
\item In the setting of algebraic number fields (\cite{Bedulev}, \cite{Topfer}).
\item If the coefficients $\ell_{i,n}$ are known to have divisors $\delta_{i,n}$ such that $\delta_{i,n}$ divides $\delta_{i,n+1}$, then (under
suitable assumptions on $\delta_{i,n}$) both $(i)$ and $(ii)$ can be improved (see  \cite{SFZu} for $(i)$,  and \cite{SFrestricted} for a detailed study of $(ii)$ with only one number, namely $r=1$). 
\item If the linear forms $\ell_{0,n} X_0 + \ell_{1,n} X_1 + \ldots+ \ell_{r,n} X_r$ are small at several points $(\xi_0^{(i)}, \ldots,  \xi_r^{(i)})$:  see \cite{SFnestsev}.
\item A lower bound can be obtained \cite{Nesterenkocritere} for the distance of the point \linebreak \mbox{$(1:\xi_1:\ldots:\xi_r)\in {\mathbb P}^r(\R)$} to any rational subspace of dimension less than $1 - \frac{\log \alpha}{\log \beta}$. This bound, proved by induction on the dimension, is Nesterenko's original approach to deduce $(i)$ from $(ii)$ (see \cite{eddzero} for a proof of this deduction based on Dirichlet's pigeonhole principle).
\end{itemize}

In all these generalizations, Proposition \ref{prop1} enables one to replace the asymptotic behaviour $\alpha^{n+o(n)}$ with the more general oscillating formula \eqref{eqcrit}.  The proof is the same as for Theorem~\ref{thcrit} (see \S \ref{subsecdemcrit}), namely one applies the result to a subsequence given by Proposition \ref{prop1}. In the case \cite{SFnestsev} where the linear forms $\ell_{0,n} X_0 + \ell_{1,n} X_1 + \ldots+ \ell_{r,n} X_r$ are small at several points, the full generality of Proposition \ref{prop1} is needed (and not only the case $N=1$).

\bigskip

At last, it would be interesting to apply Proposition \ref{prop1} to other results in Diophantine approximation. For instance, an analogous question about oscillating linear forms was asked in \cite{edd} (\S 6), and answered by Adamczewski \cite{borisedd} using ideas similar to the ones used here (but without Kronecker-Weyl's theorem); see however \cite{eddzero} for a complete answer to all questions asked in \cite{edd}.

\end{document}